\documentclass[letterpaper,10pt]{article}
\usepackage{amsmath}
\usepackage{mathtools}

\DeclareMathOperator{\tsr}{tr_{2}}
\DeclareMathOperator{\tcr}{tr_{3}}
\DeclareMathOperator{\tnr}{tr_{n}}
\usepackage{amsfonts}
\usepackage{amssymb}
\usepackage[none]{hyphenat}
\usepackage[margin=1in]{geometry}
\usepackage{fancyhdr}
\usepackage{graphicx}
\usepackage{layouts}
\graphicspath{{Images}}
\allowdisplaybreaks

\makeatletter
\newcommand\invisiblesection[1]{%
  \refstepcounter{section}%
  \addcontentsline{toc}{section}{\protect\numberline{\thesection}#1}%
  \sectionmark{#1}\phantom{}}
\makeatother

\usepackage[titles]{tocloft}

\begin{document}

\begin{titlepage}
\begin{center}
\vspace{1cm}
\Large{\textbf{Extension of Algebraic Solutions Using the Lambert W Function}}\\
\Large{\text{Sidney Edwards}}\\
\section*{Abstract}
The Lambert W function has utility for solving various exponential and logarithmic equations arranged in the form of $f(x)e^{g(x)}$. This document presents a variety of categorized inversion formulas and identities making use of the W function and tetration. Related techniques are then used to derive polar forms of exponential and related functions.
\tableofcontents
\vfill
\line(1,0){400}\\
\small{sedwards7@wisc.edu}
\end{center}
\end{titlepage}
\section{Introduction}
This document presents a series of algebraic applications of the Lambert W function, specifically in using it for deriving new solutions to various exponential and logarithmic equations. Section 2, beginning on page \pageref{Solutions} summarizes these solutions with proof of derivations starting on page \pageref{Begin Proofs}.
\\
\\
Although the W function has been defined since the time of Euler, it has only recently been applied to common sciences and mathematics. It is used in models for mechanics \cite{Cheillini}, electrical engineering \cite{OnDiode}, quantum physics \cite{step-potential} and computational science \cite{Areduction}.
\\
For example, consider an equation used in chemical engineering for modeling the time evolution of film thickness of particles in a reaction, with $D(t)$ representing film thickness over time, $a$ and $b$ representing reaction constants and W(t) is the Lambert W function \cite{Analytical}:
$\{ ( t , D(t)) \in \mathbb{R}^2 \}$
\begin{align*}
D(t) &= \frac{b}{a} \left[ 1+W(-e^{-1- a^{2}t/b^{2}}) \right]. && \\
\end{align*}
To understand how to arrive at an exact solution for time, one must know the W function's definition and properties in the following sections. The exact solution is shown on page \pageref{Exa2}.
\subsection{Overview Of The W Function}
The Lambert W function (or ``W function") denoted as $W(z)$ (for any $z \in \mathbb{C}$), is defined as the inverse of $ze^{z}$ \cite{OnThe} and satisfies the relationship $$z=W(z)e^{W(z)}.$$ 
It was originally introduced by Johann Heinrich Lambert when considering functions containing combinations of polynomials and exponents. Euler later utilized it for his famous power tower function and defined its relation to tetration \cite{TheFractal}. Although it is analytically defined in the complex plane, this document restricts the W function and consequent derivations to real-valued functions.
\begin{figure}[h]
\center
\includegraphics[scale=0.7]{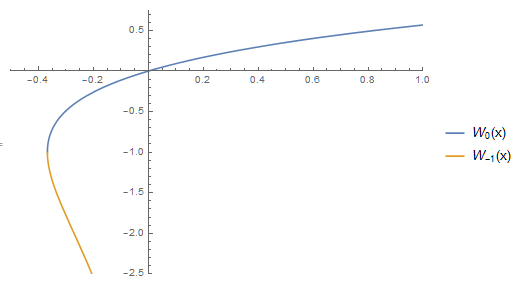}
\caption{The W function of the real variable $x$ is shown below as multivalued for $x \in {\Bbb R} $ over $[- \frac{1}{e}, 0)$ and is comprised of two branches \cite{BranchDifferences}. The principal branch, denoted $W_{0}(x)$, is the Lambert W function on the interval $[- \frac{1}{e}, \infty)$ with $W(x) \geq -1$, while the second branch, $W_{-1}(x)$, lies on the interval $[ - \frac{1}{e}, 0)$ with $W(x) \leq -1$ \cite{OnThe}.}
\end{figure}
\newpage
\subsection{Tetration}
The Lambert W function was shown to be related to tetration by Leonard Euler in his work on iterated exponentials \cite{ThePrinceton}. Tetration is a term coined for the next iterative operator beyond exponentiation or otherwise called iterated exponentiation \cite{Transfinite}. A notation similar to an exponent can be used to represent this operation:
$$
^{n}a= \underbrace{a^{a^{a^{...}}}}_\text{n times}
$$
with $a$ being the base and $n$ and being the height of tetration \cite{Areduction}. This can more simply be suggested as ``$a$ to the height of $n$". \\ \\
Tetration has two different inverse operations, one being an iterated logarithm \cite{Areduction} and other other being what is only colloquially called the super root or tetra root. The tetra root is the inverse operator that returns the base for a given $n$th height. For instance, in the case of $x^x$, or $^{2}x$, the tetra square root would be used to return $x$, $$ \tsr(x^{x}) = x.$$ 
This relates to the W function which can be used to define the tetra square root, derived from finding the inverse of $f(x)=x^x$:
\begin{align*}
y &= x^x && \\
\ln(y))  &= x \ln(x) && \\
  &= e^{ \ln(x)} \ln(x) && \\
W( \ln(y)  &= \ln(x) \\
e^{W( \ln(y)}  &= x. \\
\end{align*}
This shows that the inverse of $x^x$ is $e^{W( \ln(f(x))}$ explicitly notates the tetra square root. For real numbers, this function's principal branch is defined over the interval $[e^{ - \frac{1}{e}}, \infty)$. 
\\ \\
In working with intricate exponents, it is useful to condense the notation for certain operations. One of these is the $n$th tetra root, which can be denoted generally as $ \tnr(x)$. For clarification, $t$ represents tetra, $r$ as root and $n$ the height of tetration it inverts. \\ 
\subsection{Identities of the W function}
Starting from $$z=W(z)e^{W(z)}, \ z \in \mathbb{C} \ [1]$$ the following identities can be derived:
\\
\begin{flalign*} 
\text{W Exponential Identity} & \ \ \ \ \ \ \ \ \ \ \ \ \ \ \ \ \ \ \ \ \ \ z = W(z)e^W(z) \rightarrow & e^{W(z)} = \frac{z}{W(z)} 
\end{flalign*}
\begin{flalign*}
\text{W Log-Difference Identity} & \ \ \ \ \ \ \ \ \ \ \ \ \ \ \ \ \ \ \ln(z) = \ln \left(W(z)e^{W(z)} \right) \rightarrow & \ln(W(z) = \ln(z)-W(z)
\end{flalign*}
\begin{flalign*}
\text{W Product Identity} & \ \ \ \ \ \ \ \ \ \ \ \ \ \ \ \ \ \ \ \ \ \ \ \ \ \ \ nW(z) = W \left(nW(z)e^{nW(z)} \right)  \rightarrow & nW(z) = W \left( \frac{nz^{n}}{W(z)^{n-1}} \right)
\end{flalign*}
\begin{flalign*}
\text{W Sum Identity} & & 
\end{flalign*}
\begin{flalign*}
W(x)+W(y) = W \left([W(x)+W(y)]e^{W(x)+W(y)} \right) & \rightarrow & W(x)+W(y)= W \left(xy \left[ \frac{1}{W(y)} + \frac{1}{W(x)} \right] \right) &
\\
\end{flalign*}
A change of base formula can be derived by generalizing the base of the exponent:
\newline
\begin{minipage}[t]{0.5\textwidth}
\begin{align*}
y &= xa^{x} && \\
W_{0,a}(y) &= x && \\
y  &= xb^{ \log_{b}(a)x} && \\
\end{align*}
\end{minipage}
\begin{minipage}[t]{0.5\textwidth}
\begin{align*}
\log_{b}(a)y  &= \log_{b}(a)xb^{ \log_{b}(a)x} && \\
W_{0,b}( \log_{b}(a)y) &= \log_{b}(a)x && \\
\frac{W_{0,b}( \log_{b}(a)y)}{\log_{b}(a)} &= x && \\
\end{align*}
\end{minipage}
and thus the change of base formula is $W_{0,a}(y) = \frac{W_{0,b}( \log_{b}(a)y)}{\log_{b}(a)}.$
\section{New Solutions to Equations}
\label{Solutions}
Below is the series of various exponential and logarithmic equations with their respective inverses alluded to earlier, arranged by category. Inverses of functions similar in form to $ze^{z}$ harness the W function, starting with the basic example
$$y=xe^{ax}.$$
This function is inverted by manipulating the factor into the same form as the exponent to implore the use of the W function:
\newline
\begin{minipage}[t]{0.5\textwidth}
\begin{alignat*}{2}
y  &= xe^{ax} && \\
ay  &= axe^{ax} && \\
\end{alignat*}
\end{minipage}
\begin{minipage}[t]{0.5\textwidth}
\begin{alignat*}{2}
W(ay) &= ax && \\
\frac{W(ay)}{a} &= x. && \\
\end{alignat*}
\end{minipage}
This is also true for more complicated functions containing combinations of logarithms, exponents and specific types of polynomials (See pages \pageref{Begin Proofs}-\pageref{End Proofs} for proof of derivations).  \\ \\ In the following sections, generalized real constants are represented as $a, b, c, d, f$ and $g$ though with $f(x)$ specifically representing a function of the real variable $x$ in equations \eqref{P4}, \eqref{P5} and \eqref{P11}. ``$e$" in all derivations denotes Euler's number. In application of any following solutions, note the previously mentioned domain restrictions. Only one branch of the W function will result if the argument of the inverse function is strictly greater than 0. 
\begin{align}
\intertext{ \textbf{Products of Exponentials and Polynomials}}
\label{P1} y  &= (ax+b)^{c}e^{dx+f}  & \frac{c}{d}W \left( \frac{d}{ac}y^{\frac{1}{c}}e^{\frac{db}{ac}-\frac{f}{c}} \right) - \frac{b}{a} = &x && \\
\
\label{P2} y  &= e^{aW(x)+b}x^{c} & \frac{c}{a+c}W \left(\frac{[a+c]y^{ 1/c}}{ce^{ b/c}} \right)e^{ \frac{c}{a+c}W \left( \frac{[a+c]y^{ 1/c}}{ce^{ b/c}} \right)}   = &x && \\
\
\label{P3} y  &= ax^{b}e^{cx^{d}}+f &
\left[ \frac{bW( \frac{cd(y-f)^{d/b}}{ba^{d/b}})}{cd} \right]^{ \frac{1}{d}}  = &x && \\
\
\label{P4} f(x) &= \frac{ax+b}{ce^{-dx}-f} & \
f \left( \frac{W_{0}(\frac{cde^{bd/a}xe^{dfx/a}}{a})-W_{-1}(\frac{dfxe^{dfx/a}}{a})}{d} -bd/a \right) = &x && \\
\
\label{P5} f(x) &= ax \coth(bx)-ax &
f \left( \frac{W_{0}( \frac{-b}{a}xe^{\frac{-b}{a}})-W_{-1}(\frac{-b}{a}xe^{\frac{-b}{a}})}{2b} \right) = &x.
\
\intertext{ \textbf{Products of Logarithms and Polynomials}}
\label{P6} y  &= \frac{ \ln(ax^{b})^{c}}{dx^{f}} & \
e^{W \left( \frac{-f[a^{-f/b}dy]^{1/c}}{a^{1/b}bc} \right)} = &x && \\
\
\label{P7} y  &= ax^{b}W(x)^{c} &
( \tfrac{b+c}{b})^{1- \frac{b+c}{b}}W \left( \frac{b}{b+c} \left[ \frac{y}{a} \right]^{ \frac{1}{b+c}} \right)^{1- \frac{b+c}{b}}  = &x && \\
\
\label{P11} f(x) &= \frac{ \ln(ax+b)}{cx+d} &
f \left( \frac{e^{W((d-cb/a)xe^{(d-cb/a)x})-W(- \frac{c(d-cb/a)^{2}xe^{(d-cb/a)x}}{a})}-b}{a} \right) = &x && \\ \nonumber
\
\intertext{ \textbf{Products of Logarithms and Exponentials}}
\label{P8} y  &= W(x)e^{x} &
\ln(^{2} \tcr(e^{y}))  = &x
\
\intertext{ \textbf{Products of Logarithms}}
\label{P9} y  &= \frac{\ln(x)}{W(x)} &
\left[ \frac{1-y}{W(1-y)} \right]^{ \frac{y}{y-1}}  = &x && \\
\
\label{P10} y &= \frac{aW(bx)}{x(W(bx)+1)}+c &
\frac{ \left[1- \frac{1}{W( \frac{eab}{y-c})} \right] \frac{eab}{y-c}}{eb} = &x
\
\intertext{ \textbf{Sums of Exponentials and Polynomials}}
\label{S1} y  &= ax+b+ce^{dx} &
\frac{1}{d} \ln \left( \frac{a}{dc}W \left( \frac{dc}{a}e^{ d \frac{y-b}{a}} \right) \right)  = &x
\
\intertext{ \textbf{Sums of Logarithms and Polynomials}}
\label{S2} y  &= ax+b+ce^{W(x)} &
W \left( \frac{y-b}{a}e^{ c/a} \right)- \frac{c}{a} = & W(x) && \\
\
\label{S3} y  &= a(x+b)+c \ln \left( \frac{x}{d} \right) &
\frac{W( \frac{ad}{c}e^{ \frac{y}{c}- \frac{ab}{c}})}{a}  = &x  && \\
\
\label{S4} y  &= x^{a}+ \ln(x^{b}) &
\left[ \frac{bW(e^{ay/b})}{a} \right]^{1/a}  = &x
\
\intertext{ \textbf{Sums of Logarithms}}
\label{S5} y  &= \ln(x^{a})+W(x^{b}) &
\left[ \frac{aW(\frac{a+b}{a}e^{ by/a})}{a+b} \right]^{ \frac{1}{b}}e^{\frac{aW(\frac{a+b}{a}e^{ by/a})}{b(a+b)}} = &x && \\
\
\label{S6} y  &= W(ax)-W(bx) &
\frac{y}{abe^{y}-a^{2}} \exp \left( \frac{y}{be^{y}-a} \right) = &x && \\
\label{S7} y &= W(x)+W(1/x) & \\
& & \frac{e^{y/2}}{ \sqrt{2}} \left[y^{2}e^{y}-2 + \sqrt{(y^{2}e^{y}-2)^{2}-4} \right]^{1/2} \nonumber \\ & & \cdot e^{- \sqrt{2}e^{-y/2} \left[y^{2}e^{y}-2 + \sqrt{(y^{2}e^{y}-2)^{2}-4} \right]^{-1/2}} = x   \nonumber && \\ \nonumber \\
& & \frac{e^{-y/2}}{ \sqrt{2}} \left[y^{2}e^{y}-2 - \sqrt{(y^{2}e^{y}-2)^{2}-4} \right]^{-1/2} \nonumber \\ & & \cdot e^{\frac{e^{-y/2}}{ \sqrt{2}} \left[y^{2}e^{y}-2 - \sqrt{(y^{2}e^{y}-2)^{2}-4} \right]^{-1/2}} = x  \nonumber && \\
\
\label{S8} y &= 2W(ax)-W([bx]^{2}) & \\
& & -\frac{a}{b^{2}}e^{-y} \left(\sqrt{1-\frac{b^{2}}{a^{2}}ye^{y}}-1 \right) \nonumber \\ & & \cdot e^{-\frac{a^{2}}{b^{2}}e^{-y} \left(\sqrt{1-\frac{b^{2}}{a^{2}}ye^{y}}-1 \right)} = x \nonumber && \\
\nonumber \\
& & \frac{a}{b^{2}}e^{-y} \left( \sqrt{1-\frac{b^{2}}{a^{2}}ye^{y}}+1 \right) \nonumber \\ & & \cdot e^{ \frac{a^{2}}{b^{2}}e^{-y} \left( \sqrt{1-\frac{b^{2}}{a^{2}}ye^{y}}+1 \right)} = x \nonumber && \\
\
\label{S9} y &= x \sum_{k=1}^{n}a_{k} \ln(b_{k}x) & 
\frac{ \exp \left( W \left( \frac{ \left[ \prod_{k=1}^{n}b_{k}^{a_{k}} \right]^{1/ \sum_{k=1}^{n}a_{k}}y}{ \sum_{k=1}^{n}a_{k}} \right) \right)}{\left[ \prod_{k=1}^{n}b_{k}^{a_{k}} \right]^{1/ \sum_{k=1}^{n}a_{k}}} = &x && \\ \nonumber
\
\intertext{ \textbf{Composites of Exponentials and Polynomials}}
\label{C1} y  &= (ax+b)^{c}e^{dx+f} &
\frac{c}{d}W \left( \frac{d}{ac}y^{\frac{1}{c}}e^{\frac{db}{ac}-\frac{f}{c}} \right) - \frac{b}{a}  = &x && \\
\
\label{C2} y  &= ax^{b}e^{cx^{d}}+f &
\left[ \frac{bW( \frac{cd(y-f)^{d/b}}{ba^{d/b}})}{cd} \right]^{ \frac{1}{d}}  = &x && \\
\
\label{C3} y  &= (ax^{b})^{cx^d} &
\frac{e^{W( \frac{a^{ d/b}d}{bc} \ln(y))/d}}{a^{ 1/b}} = &x &&  \\ \nonumber
\
\intertext{ \textbf{Composites of Exponentials and Logarithms}}
\label{C4} y  &= e^{aW(x)+b}x^{c} &
\frac{c}{a+c}W \left(\frac{[a+c]y^{ 1/c}}{ce^{ b/c}} \right)e^{ \frac{c}{a+c}W \left( \frac{[a+c]y^{ 1/c}}{ce^{ b/c}} \right)}   = &x && \\
\
\label{C5} y  &= a(bx^{ \ln(c)})^{ \ln(x)} &
e^{\frac{- \ln(b) \pm \sqrt{ \ln^{2}(b)-4 \ln(c)( \ln(y)- \ln(a))}}{2 \ln(c)}}  = &x && \\
\
\label{C6} y  &= ax^{b sr_{2}(x^{c})} &
e^{ \frac{1}{c \sqrt{2}}[W(2 \ln(( \frac{y}{a})^{ \frac{c}{b}})) \ln(( \frac{y}{a})^{ \frac{c}{a}}]^{ \frac{1}{2}}} = &x && \\
\
\label{C7} y  &= (ax^{b})^{cx^d} &
\frac{ \exp \left( W \left( \frac{a^{ d/b}d}{bc} \ln(y) \right)/d \right)}{a^{ 1/b}}  = &x && \\
\label{C10} y &= a \tsr(bx)^{cx} & \frac{^{2} \left( \tcr \left( \left[ \frac{y}{a} \right]^{b/c} \right) \right)}{b} = &x && \\
\intertext{ \textbf{Composites of Logarithms and Polynomials}}
\label{C8} y  &= x^{a}+ \ln(x^{b}) &
\left[ \frac{bW(e^{ \frac{ay}{b}})}{a} \right]^{1/a}  = &x && \\
\
\label{C9} y  &= \ln(x^{a})+W(x^{b}) &
\left[ \frac{aW(\frac{a+b}{a}e^{ by/a})}{a+b} \right]^{1/b} \exp \left( \frac{aW(\frac{a+b}{a}e^{ by/a})}{b(a+b)} \right) = &x && \\ \nonumber
\
\end{align}
\subsection{Examples of Applications}
Different physical and mathematical applications can be seen harnessing the W function.
\
Starting with the diode article \cite{OnDiode}, it defines the model of the current and voltage in a solar cell it solves with the W function:
\begin{align*}
I &=I_{0}(e^{ \frac{q}{nK_{b}T}(V-IR_{s})}-1)+\frac{V-IR{s}}{R_{sh}}-I_{ph}. \ && \\
\end{align*}
Here, $I$ and $V$ denote variable current and voltage respectively. To avoid confusion, $I_0$ represents the reverse saturation current and $Iph$ is the solar cell's photocurrent. Let us denote these constants and variables more generally as this allows the equation to be more easily related to the above derivations:
\begin{align}
\label{Exa1} x+f+a &= ae^{b[y-xc]}+\frac{y-xc}{d} & \\
& & I[R_{sh}+R_{s}]+R_{sh}[I_{ph}+I_{0}]- \nonumber \\ 
& & \tfrac{n K_{b}T}{q}W \left(\tfrac{qI_{0}R_{sh}}{n K_{b}T} \exp \left( \frac{ qR_{sh}[I+I_{ph}+I_{0}]}{(n K_{b}T)} \right) \right) = V \nonumber \\
\nonumber
\end{align}
Consider an equation for modeling film thickness in a chemical reaction \cite{Analytical}. The film thickness $D(t)$ is represented as a function of time $t$ with $a$ characterizing the reaction constant in reaction \textit{a} and $b$ the diffusion in reaction \textit{b} Solving for time using properties of the W function can be done as follows:
\begin{align}
\label{Exa2} D(t) &= \frac{b}{a} \left[ 1+W(-e^{-1- a^{2}t/b^{2}}) \right] & & \\
& & \tfrac{b^{2}}{a^2}[ \tfrac{a}{b}D-1)]+ \tfrac{b^2}{a^2} \ln( \tfrac{a}{b}D-1)+ \tfrac{b^{2}}{a^{2}} = t \nonumber \\ \nonumber
\end{align}
The article on Cheillini integrability shows a solution to a differential model using the W function \cite{Cheillini}. It specifies a general equation for the non-relativistic potential $V(x)$ derived  for a quadratically damped harmonic oscillator in Lagrangian mechanics:
\begin{align*}
V^{\pm}(x_{i})= \frac{ \mu}{4}e^{-\eta \pm} \left[e^{G_{ \pm}(x_{i})}G_{ \pm}(x_{i})-e^{G_{ \pm}(x_{0})}G_{ \pm}(x_{0}) \right].
\end{align*}
$ \pm$ indicates the sign of the quadratic damping term, $x_{i}$ represents the $x$ coordinate of a turning point where the velocity of the particle is zero in phase space (indicating it models consecutive turning points), and $G_{ \pm}(x_{i})$ is a substitution for the integral of the a polynomial term with two sub-cases. After making a substitution and simplifying, the article goes on to demonstrate the steps for explicitly calculating the recursive relation $x_{i+1}$ itself using the W function:
\begin{align*}
x_{i+1}=G^{-1}_{ \pm} ( W (e^{G_{ \pm}(x_{i})}G_{ \pm}(x_{i}))).
\end{align*}
\subsection{Combining Branches of the W Function}
One relatively new concept of interest is that the W function's properties can allow both of its real branches to be used simultaneously as the identities of the W function remain invariant of the specific branch by its definition. Take for instance a function related to the generating function for Bernoulli numbers that was recently shown to be invertible using branch differences of the W function:
\begin{align*}
f(x)=\frac{x}{e^{-x}-1} \ \cite{BranchDifferences}.
\end{align*}
An algebraic method of solving this was derived by Jeffrey and Jankowski by breaking the equation into separate products using branch differences \cite{BranchDifferences}. \\ Another method of derivation can be found that exemplifies how both branches can be used from a functional perspective as in derivations \eqref{P4}, \eqref{P5} and \eqref{P11}. One should consider however that, just as with branches of other functions such as the square root function, only one definition of the function may be assumed for a given equation. \\ \\
Starting with 
\begin{align*}
 f(x)=\frac{x}{e^{-x}-1},
\end{align*}
we can substitute for a branch difference to observe how it transforms the function,
\begin{align*}
f(W_{0}(x)-W_{-1}(x))=\frac{W_{0}(x)-W_{-1}(x)}{e^{-W_{0}(x)-W_{-1}(x)}-1}. 
\end{align*}
Paying close attention to the exponent containing the branches, we can observe that it can be rewritten as 
\begin{align*}
f(W_{0}(x)-W_{-1}(x))=\frac{W_{0}(x)-W_{-1}(x)}{\frac{xW_{0}(x)}{xW_{-1}(x)}-1}
\end{align*}
which allows the $x$ in the numerator and denominator to simplify. This fraction of branches in the denominator can be manipulated and added with the constant: 
\begin{align*}
f(W_{0}(x)-W_{-1}(x)) &=\frac{W_{0}(x)-W_{-1}(x)}{\frac{xW_{0}(x)}{xW_{-1}(x)}-1} \\
f(W_{0}(x)-W_{-1}(x)) &= \frac{W_{0}(x)-W_{-1}(x)}{ \frac{W_{0}(x)-W_{-1}(x)}{W_{-1}(x)}} \\
f(W_{0}(x)-W_{-1}(x)) &= W_{-1}(x) \\  
\end{align*}
Finishing the last step requires having defined the W function in terms of the second branch to make use of the substitution $xe^{x}$,
\begin{align*}
f(W_{0}(xe^{x})-W_{-1}(xe^{x})) &= W_{-1}(xe^{x}) \\  
f(W_{0}(xe^{x})-x) &= x. \\ 
\end{align*}
From this process, one can conclude Jeffrey and Jankowski's statement that the inverse function of $ \frac{x}{e^{-x}-1}$ is $W_{0}(xe^{x})-x$. This function can be further generalized to several constants. If $f(x) = \frac{ax}{e^{-bx}-c}$, then $f^{-1}(x)$ is shown as
\begin{align}
f \left( \frac{W_{0}(\frac{bx}{a}e^{ bcx/a})-W_{-1}( \frac{bcx}{a}e^{ bcx/a})}{b} \right) &= x. \\
\nonumber
\end{align}
\section{Representation in polar coordinates}
Another instance where the W function arises is in representing Cartesian functions as polar functions. Below, one can see the conversion of several functions:
\begin{align*}
y  = \ln(x), \ y=e^x, \ y=W(x), \ y=xe^{x}.
\end{align*}
Using conventional substitution, rewrite $x$ and $y$ in terms of $r$ and $\theta$ to solve for $r$:
\\
\begin{minipage}[t]{0.5\textwidth}
\begin{align*}
y &= \ln(x) \\
r \sin( \theta) &= \ln( r \cos( \theta)) \\ 
e^{r \sin(\theta)}  &= r \cos(\theta) && \\
\frac{e^{r \sin(\theta)}}{r}  &= \cos( \theta) && \\
r e^{-r \sin(\theta)}  &= \frac{1}{ \cos( \theta)} && \\
-r e^{-r \sin(\theta)}  &= \frac{-1}{ \cos( \theta)} && \\
-r \sin(\theta) e^{-r \sin(\theta)}  &= \frac{-\sin(\theta)}{ \cos( \theta)} && \\
-r \sin(\theta) &= W(- \tan(\theta) && \\
r &= \frac{W(- \tan(\theta)}{- \sin(\theta)}. && \\
\end{align*}
\begin{align*}
y  &= xe^{x} && \\
r \sin( \theta)  &= r \cos( \theta)e^{r \cos( \theta)} && \\
\sin( \theta)  &= \cos( \theta)e^{r \cos( \theta)} && \\
\frac{ \sin( \theta)}{ \cos( \theta)}  &= e^{r \cos( \theta)} && \\
\frac{ \ln( \tan( \theta))}{ \cos( \theta)}  &= r && \\
\end{align*}
\end{minipage}
\begin{minipage}[t][][b]{0.5\textwidth}
\begin{align*}
y  &= e^{x} && \\
r \sin(\theta)  &= e^{r \cos(\theta)} && \\
\frac{1}{r \sin(\theta)}  &= e^{-r \cos(\theta)} && \\
\frac{1}{ \sin(\theta)}  &= r e^{-r \cos(\theta)} && \\
-\frac{ \cos(\theta)}{ \sin(\theta)}  &= -r \cos(\theta)e^{-r \cos(\theta)} && \\
W(- \cot(\theta))  &= -r \cos(\theta) && \\
\frac{W(- \cot(\theta))}{- \cos(\theta)}  &= r  && \\
\end{align*}
\begin{align*}
y  &= W(x) && \\
r \sin(\theta)  &= W(r \cos(\theta)) && \\
r \sin(\theta)e^{r \sin(\theta)}  &= r \cos(\theta) && \\
\sin(\theta)e^{r \sin(\theta)}  &= \cos(\theta) && \\
e^{r \sin(\theta)}  &= \frac{ \cos(\theta)}{ \sin(\theta)} && \\
r \sin(\theta) &= \ln( \cot(\theta)) && \\
r &= \frac{ \ln( \cot(\theta))}{ \sin(\theta)} && \\
\end{align*}
\end{minipage}
From these representations, it is easy to see one can use polar coordinates to transform a function such as $\ln(x)$ to its Cartesian inverse 
\begin{align*}
\frac{W(- \tan( \pi - [\theta + \frac{ \pi}{2}]))}{-\sin(  \pi - [ \theta + \frac{ \pi}{2}])} &= \frac{W(- \cot( \theta))}{- \cos( \theta)}. 
\end{align*}
\subsection{Rotations in Cartesian coordinates}
Moving further with this notion, one can also use the polar form to solve for rotations of the curves created by these equations in Cartesian coordinates. \\ Take for instance the equation $y=e^{x},$ put into terms of $r$ and $ \theta$ as $r \sin( \theta) = e^{r \cos( \theta)}.$ Starting from this, a rotation of the curve by $ \frac{ \pi}{4}$ can be induced via the argument which then allows one to solve for the rotated curve in terms of Cartesian coordinates using the W function:
\begin{align*}
r \sin( \theta) &= e^{r \cos( \theta)} \rightarrow r \sin( \theta- \frac{ \pi}{4}) = e^{r \cos( \theta- \pi/4)}.
\end{align*}
In Figure 2 on page \pageref{polar figure} are the graphics for these implicit curves. The left image is of the unrotated exponential function (in polar form) while the right image shows this same curve but rotated by $ \frac{ \pi}{4}$ in polar coordinates.
\begin{figure}[h]
\center
\includegraphics[scale=.7]{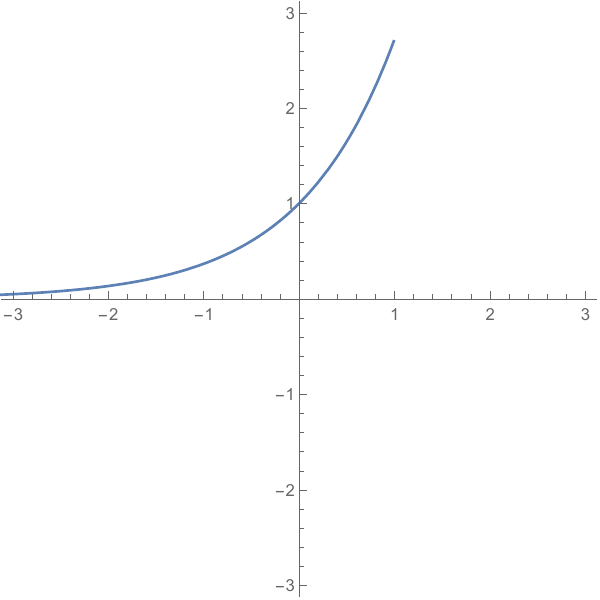}
\includegraphics[scale=.7]{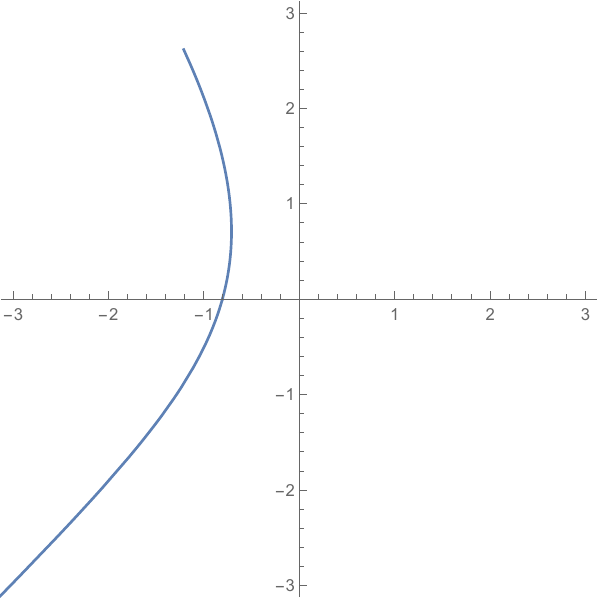}
\caption{The original exponential function $y=e^{x}$ shown in the left image is plotted in polar form in Mathematica using the preceding equations. The left image contains the same curve but rotated in polar coordinates by $ \frac{ \pi}{4}$.}
\label{polar figure} 
\end{figure}
Using conventional trigonometric identities, the right-hand side can be expanded out
\begin{align*}
r \sin( \theta- \frac{ \pi}{4}) &= e^{r \cos( \theta-\frac{ \pi}{4})} \rightarrow r( \frac{ \sin( \theta)}{ \sqrt{2}}- \frac{ \cos( \theta)}{ \sqrt{2}}) = e^{r( \frac{sin( \theta)}{ \sqrt{2}}+ \frac{ \cos( \theta)}{ \sqrt{2}})}
\end{align*}
It is convenient here to make two substitutions of $u=\frac{r \sin( \theta)}{ \sqrt{2}}$ and $v=\frac{r \cos( \theta)}{ \sqrt{2}}$, transforming the preceding equation into the form of 
\begin{align*}
u-v=e^{u+v}.
\end{align*}
We have seen a similar equation to this in \eqref{S1} that is solved using the W function. A similar technique can be applied in this specific instance to solve for $u$ and thus for its related Cartesian form $y$:
\begin{align*}
-v &= e^{u+v}-u && \\
-2v &= e^{u+v}-u-v && \\
2v &= -e^{u+v}+u+v && \\
e^{2v} &= e^{-e^{u+v}+u+v} && \\
e^{2v} &= e^{u+v}e^{-e^{u+v}} && \\
-e^{2v} &= -e^{u+v}e^{-e^{u+v}} && \\
W(-e^{2v}) &= -e^{u+v} && \\
\ln(-W(-e^{2v})) &= u+v && \\
\ln(-W(-e^{2v}))-v &= u && \\
\ln \left( -W \left(- \exp \left( 2\frac{r \cos( \theta)}{ \sqrt{2}} \right) \right) \right)-\tfrac{r \cos( \theta)}{ \sqrt{2}} &= \tfrac{r \sin( \theta)}{ \sqrt{2}} && \\
\sqrt{2} \ln(-W(- \exp( \sqrt{2}x )))-x &= y && \\
\end{align*}
The top-most graph in Figure \ref{Graph2} on page \pageref{Graph2} shows the plots in polar form. Below shows the plots in their Cartesian form. We can see on the top of the two graphs that the original $y=e^x$ is plotted while on the right the same curve is rotated by $ \frac{ \pi}{4}$ and plotted in the form derived from above.  \\ \\ This rotation can be generalized to any angle "$ \phi$" in a similar manner with the exception of angles of $\frac{ \pi n}{2}$.
\begin{figure}[h]
\label{Graph2}
\includegraphics[scale=.5]{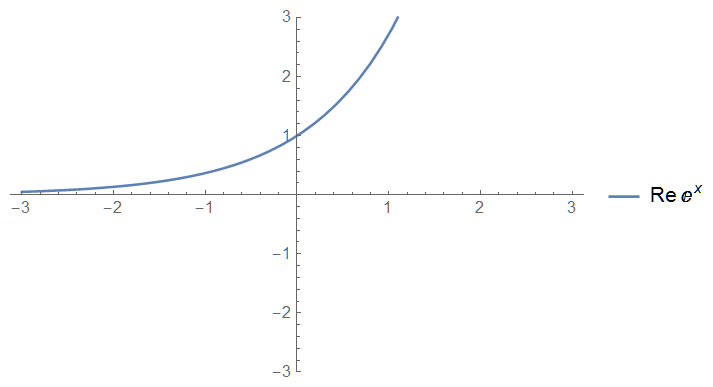}
\includegraphics[scale=.5]{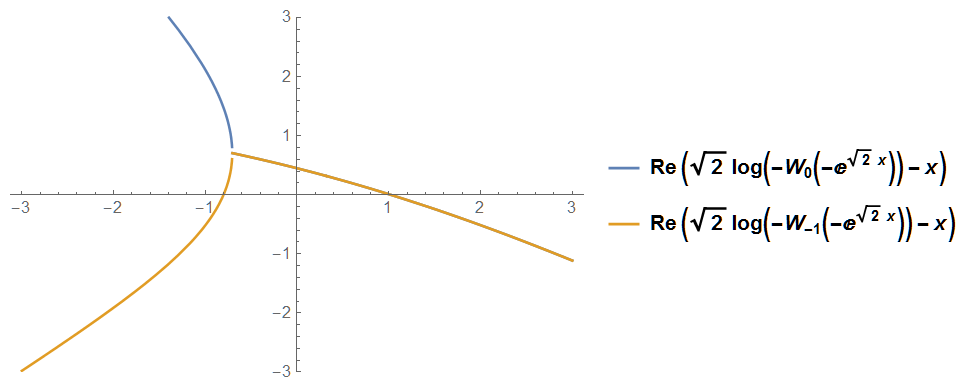}
\caption{The two images show the standard exponential function in Cartesian form in the upper graph. In the lower graph is a plot of the curve rotated by $ \frac{ \pi}{4}$ expressed in Cartesian form. The left concave curve represents where the real component of both branches of the rotated function overlap.} 
\end{figure}  
Starting with the original equation, $y=e^{x}$, it is represented in polar coordinates as
$r \sin{ \theta} = e^{r \cos( \theta)}$
which can be rotated by an arbitrary angle $ \phi$ (for $ \phi \neq \frac{ n \pi}{2}$) and solved for $y$ in Cartesian form:
\begin{align}
\label{CTP1} y &= Ae^{Bx}  && \\
\rightarrow r \sin( \theta+ \phi) &= Ae^{Br \cos( \theta+ \phi)}  \nonumber && \\
u &= r \cos(\theta),  \nonumber &&\\
a &= \cos( \phi), \nonumber &&\\
v &= r \sin(\theta), \nonumber &&\\
b &= \sin( \phi), \nonumber &&\\
\tfrac{1}{Bb} \ln( \tfrac{a}{ABb}W( \tfrac{ABb}{a} \exp( B\frac{a^2+b^2}{a}u)))- \frac{a}{b}u &= v  \nonumber && \\
\tfrac{ \csc( \phi)}{B} \ln( \tfrac{ \cot( \phi)}{AB}W(AB \tan( \phi)e^{B \sec( \phi)x}))- \cot( \phi)x &= y  \nonumber && \\
\nonumber
\end{align}
Both branches of the W function are needed to describe the real rotation of curve excluding when the argument of W is strictly greater than or equal to zero.
Almost the same exact process can be used to find the rotation for the natural logarithm function in Cartesian form as well:
\begin{align}
\label{CTP2} y &= \frac{1}{B} \ln \left( \frac{x}{A} \right) && \\
\rightarrow \sin( \theta+ \phi) &= \tfrac{1}{B} \ln \left( \frac{r \cos( \theta+ \phi)}{A} \right) \nonumber && \\
u &= r \cos(\theta), \nonumber &&\\
a &= \cos( \phi), \nonumber &&\\
v &= r \sin(\theta), \nonumber &&\\
b &= \sin( \phi), \nonumber &&\\
v &= \tfrac{1}{Ba} \ln( \tfrac{b}{ABa} W( \tfrac{ABa}{b}e^{\frac{Ba^{2}+Bb^{2}}{b}u}))- \frac{b}{a}u \nonumber && \\
y &= \tfrac{ \sec( \phi)}{B} \ln( \tfrac{ \tan( \phi)}{AB}W(AB \cot( \phi) e^{B \csc( \phi)x}))- \tan( \phi)x   \nonumber && \\ \nonumber
\end{align}
One may notice that the trigonometric functions within the rotations of both $e^{x}$ and $ \ln(x)$ have recurring singularities for certain values of $ \phi$. The Cartesian rotation equation is invalid for interval rotations of $ \phi = \frac{ \pi n}{2}$. The proper procedure for $ \phi = \frac{ \pi n}{2}$ should therefore be to derive the Cartesian rotation equation separately. In the instance of $ \phi = \frac{ \pi}{2}$, the resulting back-conversion to Cartesian form is $y=- \ln(x)$.
\section{Derivations}
\label{Begin Proofs}
The following pages contain the algebraic proofs of derivations in previous sections of this document. The first section of derivations corresponds to section 2 on page \pageref{Solutions}, arranged in the same order but without repeating identical derivations. Later derivations are proofs of the example applications in section 2 and polar representations in section 3. \\ \\
When working with such formulas, it is important to care for the domain restrictions of the W function. Since $W(x)$ is not bijective, $W(xe^x)=x$ cannot always be assumed \cite{OnThe}. The principal branch $W_{0}(x)$ is real only over $x \in [-1/e, \infty)$ and since $W(-1e^{-1})=-1$,  it is necessary to specify the constraint that $g(x) \geq -1$ for $W_{0}(g(x)e^{g(x)})=g(x)$ to hold. Similarly, it is also necessary to that $g(x) \leq -1$ for $W_{-1}(g(x)e^{g(x)})=g(x)$ to hold. \\ \\
In derivation \eqref{P1}, the inverse function of $y$ is proven as $\frac{W \left(y^{ \frac{a}{cd}}e^{- \frac{af}{d}+b} \right)}{a}.$ Since the domain of the principal branch is real for $x \in [-1/e, \infty)$,  the argument
$y^{ \frac{a}{cd}}e^{- \frac{af}{d}+b}$ must also be greater than -1/e.
\begin{alignat*}{2}
\frac{y^{ \frac{a}{cd}}}{e^{ \frac{af}{d}-b}} \geq -1/e,  \ \ \ \ \ & \ \ \ \ \ y^{ \frac{a}{cd}} \geq -e^{ \frac{af}{d}-b-1}, & \ \ \ \ \ \ \ \ \ \ y & \geq (-1)^{\frac{cd}{a}}e^{ \frac{acf-bcd-cd}{a}}.
\end{alignat*}
This reveals that one must require $\left\{ \frac{cd}{a} \in \mathbb{Z} \ | \ ( \exists z) \ \frac{cd}{a}=z \right\}$ to avoid complex values in the argument. One must also constrain the exponentiated term $\frac{acf-bcd-cb}{a}$ such that $ (\forall \frac{cd}{a})=2z+1, \frac{acf-bcd-cd}{a} \geq -1$. In the case that $ \frac{cd}{a}=2z$, the exponential has a positive coefficient and thus the exponentiated term can be taken as any real number that satisfies the first constraint. Similarly for $W_{-1}$, $\frac{cd}{a}=2z+1$ for the argument of $W_{-1}$ to remain in $[-1/e,0)$.
\newpage
\begin{minipage}[t]{0.5\textwidth} 
\begin{flalign*}
\intertext{Derivation \eqref{P1}, \eqref{C1}}
y  &= (ax+b)^{c}e^{dx+f} && \\
y^{1/c} &=(ax+b)e^{\frac{dx}{c}+\frac{f}{c}} && \\
 y^{1/c}e^{-f/c} &= (ax+b)e^{dx/c} && \\
\frac{d}{ac}y^{1/c}e^{-f/c} &= \frac{d}{ac}(ax+b)e^{dx/c} && \\
&= \left( \frac{dx}{c} + \frac{db}{ac} \right)e^{dx/c} && \\
\frac{d}{ac}y^{1/c}e^{\frac{db}{ac}-\frac{f}{c}} &= \left( \frac{dx}{c}+\frac{db}{ac} \right) e^{\frac{dx}{c}+\frac{db}{ac}} && \\
W \left( \frac{d}{ac}y^{1/c}e^{\frac{db}{ac}-\frac{f}{c}} \right) &= \frac{dx}{c} + \frac{db}{ac}  && \\
W \left( \frac{d}{ac}y^{1/c}e^{\frac{db}{ac}-\frac{f}{c}} \right) - \frac{db}{ac} &= \frac{dx}{c} && \\
\frac{c}{d}W \left( \frac{d}{ac}y^{1/c}e^{\frac{db}{ac}-\frac{f}{c}} \right) - \frac{b}{a} &= x
\end{flalign*}
\end{minipage}
\begin{minipage}[t]{0.5\textwidth} 
\begin{flalign*}
\intertext{Derivation \eqref{P2}, \eqref{C4}}
y  &= e^{aW(x)+b}x^{c} && \\
  &= e^{aW(x)+b} \left(W(x)e^{W(x)} \right)^{c} && \\
  &= e^{aW(x)+b}e^{cW(x)}W(x)^{c} && \\
\frac{y}{e^{b}}  &= e^{aW(x)}e^{cW(x)}W(x)^{c} && \\
  &= e^{aW(x)+cW(x)}W(x)^{c} && \\
\frac{y^{ \frac{1}{c}}}{e^{ \frac{b}{c}}}  &= e^{ \frac{[a+c]}{c}W(x)}W(x) && \\
\frac{[a+c]y^{ 1/c}}{ce^{ b/c}}  &= e^{ \frac{[a+c]}{c}W(x)}W(x) \tfrac{a+c}{c} && \\
W \left( \frac{[a+c]y^{ 1/c}}{ce^{ b/c}} \right)  &= W(x) \tfrac{a+c}{c} && \\
\frac{cW \left( \frac{[a+c]y^{ 1/c}}{ce^{ \frac{b}{c}}} \right)}{a+c}  &= W(x) && \\
\frac{cW \left(\frac{[a+c]y^{ 1/c}}{ce^{ \frac{b}{c}}} \right)}{a+c}e^{\frac{cW \left( \frac{[a+c]y^{ 1/c}}{ce^{ \frac{b}{c}}} \right)}{a+c}}   &= x && \\
\end{flalign*}
\newpage
\end{minipage}
\begin{flalign*}
\intertext{Derivation \eqref{P3}, \eqref{C2}}
y  &= ax^{b}e^{cx^{d}}+f && \\
y-f  &= ax^{b}e^{cx^{d}} && \\
\frac{y-f}{a}  &= x^{b}e^{cx^{d}} && \\
\frac{(y-f)^{1/b}}{a^{1/b}}  &= xe^{ \frac{cx^{d}}{b}} && \\
\frac{(y-f)^{d/b}}{a^{ \frac{d}{b}}}  &= x^{d}e^{ \frac{cdx^{d}}{b}} && \\
\frac{cd(y-f)^{d/b}}{ba^{ \frac{d}{b}}}  &= \frac{cdx^{d}}{b}e^{ \frac{cdx^{d}}{b}} && \\
W \left( \frac{cd(y-f)^{d/b}}{ba^{ \frac{d}{b}}} \right)  &= \frac{cdx^{d}}{b} && \\
\frac{bW \left( \frac{cd(y-f)^{d/b}}{ba^{ \frac{d}{b}}} \right)}{cd}  &= x^d && \\
\left[ \frac{bW \left( \frac{cd(y-f)^{ d/b}}{ba^{d/b}} \right)}{cd} \right]^{ 1/d}  &= x && \\
\end{flalign*}
\newpage
\begin{align*}
\intertext{Derivation \eqref{P4}}
f(x) &= \frac{ax+b}{ce^{-dx}-f}. && \\
f \left( \frac{W_{0}(x)-W_{-1}(x)}{d} \right) &= \frac{ \tfrac{a}{d}[W_{0}(x)-W_{-1}(x)]+b}{ce^{-[W_{0}(x)-W_{-1}(x)]}-f}. 
&& \\
f \left( \frac{W_{0}(x)-W_{-1}(x)}{d} -bd/a \right) &= \frac{ \tfrac{a}{d}[W_{0}(x)-W_{-1}(x)-bd/a]+b}{ce^{-[W_{0}(x)-W_{-1}(x)]}-f}. 
&& \\
&= \frac{ \tfrac{a}{d}[W_{0}(x)-W_{-1}(x)]}{ce^{-[W_{0}(x)-W_{-1}(x)]-bd/a}-f}. 
&& \\
&= \frac{ \tfrac{a}{d}[W_{0}(x)-W_{-1}(x)]}{ce^{bd/a}\frac{W_{0}(x)}{W_{-1}(x)}-f}. 
&& \\
f \left( \frac{W_{0}(x)-W_{-1}(fx)}{d} -bd/a \right) &= \frac{ \tfrac{a}{d}[W_{0}(x)-W_{-1}(fx)]}{cfe^{bd/a}\frac{W_{0}(x)}{W_{-1}(fx)}-f}. 
&& \\
&= \frac{ \tfrac{a}{d}[W_{0}(x)-W_{-1}(fx)]}{f[ce^{bd/a}\frac{W_{0}(x)}{W_{-1}(fx)}-1]} 
&& \\
f \left( \frac{W_{0}(ce^{bd/a}x)-W_{-1}(fx)}{d} -bd/a \right) &= \frac{ \tfrac{a}{d}[W_{0}(ce^{bd/a}x)-W_{-1}(fx)]}{f[ \frac{W_{0}(ce^{bd/a}x)}{W_{-1}(fx)}-1]} 
&& \\
&= \frac{ \tfrac{a}{d}[W_{0}(ce^{bd/a}x)-W_{-1}(fx)]}{f \left[ \frac{W_{0}(ce^{bd/a}x)-W_{-1}(fx)}{W_{-1}(fx)} \right]} 
&& \\
f \left( \frac{W_{0}(ce^{bd/a}xe^{fx})-W_{-1}(fxe^{fx})}{d} -bd/a \right) &= \frac{a}{d}x
&& \\
f \left( \frac{W_{0}(\frac{cde^{bd/a}xe^{dfx/a}}{a})-W_{-1}(\frac{dfxe^{dfx/a}}{a})}{d} -bd/a \right) &= x
&&\\
\end{align*}
\newpage
\begin{align*}
\intertext{Derivation \eqref{P5}}
f(x) &= ax \coth(bx)+ax && \\
&= ax \left( \coth(bx)+1 \right) && \\
f( \frac{x}{2b}) &= \frac{ax}{2b} \left( \coth(x/2)+1 \right) && \\
&= \frac{ax}{2b}( \frac{e^{x/2}+e^{-x/2}}{e^{x/2}-e^{-x/2}}+1) && \\
&= \frac{ax}{2b} \left( \frac{e^{x/2}+e^{-x/2}}{e^{x/2}-e^{-x/2}}+ \frac{e^{x/2}-e^{-x/2}}{e^{x/2}-e^{-x/2}} \right) && \\
&= \frac{ax}{2b} \left( \frac{2e^{x/2}}{e^{x/2}-e^{-x/2}} \right) && \\
&= \frac{ax}{b} \left( \frac{e^{x/2}}{e^{x/2}-e^{-x/2}} \right) && \\
f \left( \frac{W_{0}(x)-W_{-1}(x)}{2b} \right) &= \frac{a(W_{0}(x)-W_{-1}(x))}{b} \left( \frac{e^{(W_{0}(x)-W_{-1}(x))/2}}{(e^{W_{0}(x)-W_{-1}(x))/2}-e^{-(W_{0}(x)+W_{-1}(x))/2}} \right) && \\
&= \frac{a(W_{0}(x)-W_{-1}(x))}{b} \left( \frac{ \sqrt{ \frac{W_{-1}(x)}{W_{0}(x)}}}{ \sqrt{ \frac{W_{-1}(x)}{W_{0}(x)}}-\sqrt{ \frac{W_{0}(x)}{W_{-1}(x)}}} \right) && \\
&= \frac{a(W_{0}(x)-W_{-1}(x))}{b} \left( \frac{ \sqrt{ \frac{W_{-1}(x)}{W_{0}(x)}}}{ \frac{W_{-1}(x)-W_{0}(x)}{ \sqrt{W_{0}(x)W_{-1}(x)}}} \right) && \\
&= \frac{-a}{b} \left( \frac{ \sqrt{ \frac{W_{-1}(x)}{W_{0}(x)}}}{ \frac{1}{ \sqrt{W_{0}(x)W_{-1}(x)}}} \right) && \\
&= \frac{-a}{b} \sqrt{ \frac{W_{0}(x)W_{-1}(x)^{2}}{W_{0}}} && \\
&= \frac{-a}{b}W_{-1}(x) && \\
f \left( \frac{W_{0}( \frac{-b}{a}xe^{\frac{-b}{a}})-W_{-1}(\frac{-b}{a}xe^{\frac{-b}{a}})}{2b} \right) &= \frac{-a}{b}W_{-1}\left( \frac{-b}{a}xe^{\frac{-b}{a}} \right) && \\
&= \frac{-a}{b} \frac{(-b)}{a}x && \\
&= x && \\
\end{align*}
\newpage
\begin{minipage}[t]{0.2\textwidth}
\begin{flalign*}
\intertext{Derivation \eqref{P6}}
y  &= \frac{ \ln(ax^{b})^{c}}{dx^{f}} && \\
dy  &= \frac{ \ln(ax^{b})^{c}}{x^{f}} && \\
(dy)^{ \frac{1}{c}} &= \frac{ \ln(ax^{b})}{x^{ \frac{f}{c}}} && \\
(dy)^{ \frac{1}{c}} &= \frac{ \ln(ax^{b})}{e^{ \ln(x) \frac{f}{c}}} && \\
(dy)^{ \frac{1}{c}} &= \ln(ax^{b})e^{\frac{-f\ln(x)}{c}} && \\
(dy)^{ \frac{1}{c}} &= b \ln(a^{ \frac{1}{b}}x)e^{\frac{-f\ln(x)}{c}} && \\
\frac{e^{ \frac{-f\ln(a^{ \frac{1}{b}})}{c}}(dy)^{ \frac{1}{c}}}{b} &= \ln(a^{ \frac{1}{b}}x)e^{\frac{-f( \ln(x)+ \ln(a^{ \frac{1}{b}})}{c}} && \\
\frac{(a^{- \frac{f}{b}}dy)^{ \frac{1}{c}}}{b} &= \ln(a^{ \frac{1}{b}}x)e^{\frac{-f\ln(a^{ \frac{1}{b}}x)}{c} } && \\
\tfrac{-f(a^{- \frac{f}{b}}dy)^{ \frac{1}{c}}}{bc} &= \tfrac{-f}{c}\ln(a^{ \frac{1}{b}}x)e^{\frac{-f \ln(a^{ \frac{1}{b}}x)}{c}} && \\
W(\frac{-f(a^{- \frac{f}{b}}dy)^{ \frac{1}{c}}}{bc}) &= \ln(a^{ \frac{1}{b}}x) && \\
e^{W(\frac{-f(a^{- \frac{f}{b}}dy)^{ \frac{1}{c}}}{bc})} &= a^{ \frac{1}{b}}x && \\
e^{W \left( \frac{-f(a^{- \frac{f}{b}}dy)^{ \frac{1}{c}}}{a^{ \frac{1}{b}}bc} \right)} &= x && \\
\end{flalign*}
\end{minipage}
\begin{minipage}[t]{0.2\textwidth}
\begin{flalign*}
\intertext{Derivation \eqref{P7}}
y  &= ax^{b}W(cx^d)^{f} 
&& \\
\frac{y}{a} &= x^{b}W(cx^d)^{f}
&& \\ 
\left[ \frac{y}{a} \right]^{1/b} &= xW(cx^d)^{f/b}
&& \\
\left[ \frac{y}{a} \right]^{d/b} &= x^{d}W(cx^d)^{df/b}
&& \\
c \left[ \frac{y}{a} \right]^{d/b} &= cx^{d}W(cx^d)^{df/b}
&& \\
&= W(cx^{d})e^{W(cx^{d})}W(cx^d)^{\frac{df}{b}}
&& \\
&= e^{W(cx^{d})}W(cx^d)^{df/b+1}
&& \\
&= e^{W(cx^{d})}W(cx^d)^{(df+b)/b}
&& \\
c \left[ \frac{y}{a} \right]^{d/(df+b)} &= \exp \left( \frac{bW(cx^{d})}{df+b} \right)W(cx^d)
&& \\
\frac{bc}{df+b} \left[ \frac{y}{a} \right]^{d/(df+b)} &= \frac{bW(cx^d)}{df+b} \exp \left( \frac{bW(cx^{d})}{df+b} \right)
&& \\
W \left( \frac{bc}{df+b} \left[ \frac{y}{a} \right]^{d/(df+b)} \right) &= \frac{bW(cx^d)}{df+b}
&& \\
\frac{df+b}{b}W \left( \frac{bc}{df+b} \left[ \frac{y}{a} \right]^{d/(df+b)} \right) &= W(cx^d)
&& \\
\tfrac{df+b}{b}W \left( \tfrac{bc}{df+b} \left[ \frac{y}{a} \right]^{d/(df+b)} \right) & \\ \boldsymbol{\cdot} \exp \left( \tfrac{df+b}{b}W \left( \tfrac{bc}{df+b} \left[ \frac{y}{a} \right]^{d/(df+b)} \right) \right) &= cx^d
&& \\
\tfrac{df+b}{bc}W \left( \tfrac{bc}{df+b} \left[ \frac{y}{a} \right]^{d/(df+b)} \right) & \\ \boldsymbol{\cdot} \exp \left( \tfrac{df+b}{b}W \left( \tfrac{bc}{df+b} \left[ \frac{y}{a} \right]^{d/(df+b)} \right) \right) &= x^d
&& \\
\left[ \tfrac{df+b}{bc}W \left( \tfrac{bc}{df+b} \left[ \frac{y}{a} \right]^{\frac{d}{df+b}} \right) \right]^{\frac{1}{d}} & \\ \boldsymbol{\cdot} \left[ \exp \left( \tfrac{df+b}{b}W \left( \tfrac{bc}{df+b} \left[ \frac{y}{a} \right]^{\frac{d}{df+b}} \right) \right) \right]^{\frac{1}{d}} &= x
&& \\
\end{flalign*}
\end{minipage}
\newpage
\begin{flalign*}
\intertext{Derivation \eqref{P11}}
f(x) &= \frac{ \ln(ax+b)}{cx+d} 
&& \\
f( \frac{x-b}{a}) &= \frac{ \ln(x)}{c \left( \frac{x-b}{a} \right)+d} 
&& \\
&= \frac{ \ln(x)}{ \frac{c}{a}x+(d-\frac{cb}{a})} 
&& \\
&= \frac{ \ln(x)}{(d- \frac{cb}{a}) \left( \frac{c}{a(d-cb/a)}x+1 \right)} 
&& \\
f( \frac{e^{x}-b}{a}) &= \frac{ x}{(d- \frac{cb}{a}) \left( \frac{ce^{x}}{a(d-cb/a)}+1 \right)} 
&& \\
f( \frac{e^{W(ax)-W(cx)}-b}{a}) &= \frac{W(ax)-W(cx)}{(d- \frac{cb}{a}) \left( \frac{ce^{W(ax)-W(cx)}}{a(d-cb/a)}+1 \right)} 
&& \\
&= \frac{W(ax)-W(cx)}{(d- \frac{cb}{a}) \left( \frac{W(cx)}{(d-cb/a)W(ax)}+1 \right)} 
&& \\
f \left( \frac{e^{W( \frac{ax}{d-cb/a})-W(cx)}-b}{a} \right) &= \frac{W( \frac{ax}{d-cb/a})-W(cx)}{(d- \frac{cb}{a}) \left( \frac{W(cx)}{W( \frac{ax}{d-cb/a})}+1 \right)} 
&& \\
f \left( \frac{e^{W( \frac{ax}{d-cb/a})-W(-cx)}-b}{a} \right) &= \frac{W( \frac{ax}{d-cb/a})-W(-cx)}{(d- \frac{cb}{a}) \left( \frac{-W(-cx)+W( \frac{ax}{d-cb/a})}{W( \frac{ax}{d-cb/a})} \right)} 
&& \\
&= \frac{W \left( \frac{ax}{d-cb/a} \right)}{d- \frac{cb}{a}}
&& \\\
f \left( \frac{e^{W((d-cb/a)xe^{(d-cb/a)x})-W(- \frac{c(d-cb/a)^{2}xe^{(d-cb/a)x}}{a})}-b}{a} \right) &= x
\end{flalign*}
\newpage
\begin{minipage}[t]{0.5\textwidth}
\begin{flalign*}
\intertext{Derivation \eqref{P8}}
y  &= W(x)e^{x} && \\
y  &= W(x)e^{W(x)e^{W(x)}} && \\
e^{y}  &= e^{W(x)e^{W(x)e^{W(x)}}} && \\
\tcr(e^{y})  &= e^{W(x)} && \\
\ln(\tcr(e^{y}))  &= W(x) && \\
\ln(^{2} \tcr(e^{y}))  &= x && \\
\end{flalign*}
\begin{flalign*}
\intertext{Derivation \eqref{P9}}
y  &= \frac{\ln(x)}{W(x)} && \\
u  &= W(x) && \\
y  &= \frac{ \ln(ue^{u})}{u} && \\
y  &= \frac{ \ln(u)+u}{u} && \\
y  &= \frac{ \ln(u)}{u}+1 && \\\\
y-1  &= \frac{ \ln(u)}{u} && \\
v  &= \ln(u) &&\\
y-1  &= ve^{-v} && \\
1-y  &= -ve^{-v} && \\
W(1-y)  &= -v && \\
-W(1-y)  &= v && \\
-W(1-y)  &= \ln(u) && \\
-W(1-y)  &= \ln(W(x)) && \\
e^{-W(1-y)}  &= W(x) && \\
e^{-W(1-y)}e^{e^{-W(1-y)}}  &= x && \\
\left( \frac{1-y}{W(1-y)} \right)^{ \frac{y}{y-1}}  &= x && \\
\end{flalign*}
\end{minipage}
\begin{minipage}[t]{0.5\textwidth}
\begin{align*}
\intertext{Derivation \eqref{P10}}
y &= \frac{aW(bx)}{x[W(bx)+1]}+c && \\
y-c &= \frac{abe^{-W(bx)}}{[W(bx)+1]} && \\
\frac{1}{y-c} &= \frac{e^{W(bx)}[W(bx)+1]}{ab} && \\
\frac{e}{y-c} &= \frac{e^{W(bx)+1}[W(bx)+1]}{ab} && \\
\frac{eab}{y-c} &= e^{W(bx)+1}[W(bx)+1] && \\
W( \tfrac{eab}{y-c}) &= W(bx)+1 && \\
W( \tfrac{eab}{y-c})-1 &= W(bx) && \\
[W \left( \tfrac{eab}{y-c} \right)-1]e^{[W( \frac{eab}{y-c})-1]} &= bx && \\
\frac{[W( \frac{eab}{y-c})-1]e^{[W( \frac{eab}{y-c})-1]}}{b} &= x && \\
\frac{[W( \frac{eab}{y-c})-1]e^{W( \frac{eab}{y-c})}}{eb} &= x && \\
\frac{[W( \frac{eab}{y-c})-1]\frac{ \frac{eab}{y-c}}{W( \frac{eab}{y-c})}}{eb} &= x && \\
\frac{[1- \frac{1}{W( \frac{eab}{y-c})}] \frac{eab}{y-c}}{eb} &= x && \\
\end{align*}
\end{minipage}
\newpage
\begin{minipage}[t]{0.5\textwidth}
\begin{flalign*}
\intertext{Derivation \eqref{S1}}
y  &= ax+b+ce^{dx} &&\\
u  &= e^{x} &&\\
y  &= a \ln(u)+b+cu^{d} &&\\
y-b  &= a \ln(u)+cu^{d} &&\\
e^{y-b}  &= e^{a \ln(u)+cu^{d}} &&\\
e^{y-b}  &= u^{a}e^{cu^{d}} &&\\
e^{ \frac{y-b}{a}}  &= ue^{ \frac{c}{a}u^{d}} &&\\
e^{ d \frac{y-b}{a}}  &= u^{d}e^{ \frac{cd}{a}u^{d}} &&\\
\frac{cd}{a}e^{ d \frac{y-b}{a}}  &= \frac{cd}{a}u^{d}e^{ \frac{cd}{a}u^{d}} &&\\
W(\frac{cd}{a}e^{ d \frac{y-b}{a}})  &= \frac{dc}{a}u^{d} &&\\
\frac{a}{cd}W(\frac{cd}{a}e^{ d \frac{y-b}{a}})  &= u^{d} &&\\
(\frac{a}{cd}W(\frac{cd}{a}e^{ d \frac{y-b}{a}}))^{ \frac{1}{d}}  &= u &&\\
(\frac{a}{cd}W(\frac{cd}{a}e^{ d \frac{y-b}{a}}))^{ \frac{1}{d}}  &= e^{x} &&\\
\frac{1}{d} \ln(\frac{a}{cd}W(\frac{cd}{a}e^{ d \frac{y-b}{a}}))  &= x &&\\
\end{flalign*}
\begin{flalign*}
\intertext{Derivation \eqref{S2}}
y  &= ax+b+ce^{W(x)} &&\\
y-b  &= ax+ce^{W(x)} &&\\
y-b  &= aW(x)e^{W(x)}+ce^{W(x)} &&\\
y-b  &= e^{W(x)}(aW(x)+c) &&\\
y-b  &= e^{W(x)}(W(x)+\frac{c}{a})a &&\\
\frac{y-b}{a}  &= e^{W(x)}(W(x)+\frac{c}{a}) &&\\
\frac{y-b}{a}e^{ \frac{c}{a}}  &= e^{ \frac{c}{a}}e^{W(x)}(W(x)+\frac{c}{a}) &&\\
\frac{y-b}{a}e^{ \frac{c}{a}}  &= e^{W(x)+\frac{c}{a}}(W(x)+\frac{c}{a}) &&\\
W( \frac{y-b}{a}e^{ \frac{c}{a}})  &= W(x)+\frac{c}{a} &&\\
W( \frac{y-b}{a}e^{ \frac{c}{a}})- \frac{c}{a}  &= W(x) &&\\
W( \frac{y-b}{a}e^{ \frac{c}{a}})- \frac{c}{a}  &= W(x) &&\\
\end{flalign*}
\end{minipage}
\begin{minipage}[t]{0.5\textwidth}
\begin{flalign*}
\intertext{Derivation \eqref{S3}}
y  &= a(x+b)+c \ln( \tfrac{x}{d}) && \\
e^{y}  &= \frac{x^{c}}{d^{c}}e^{a(x+b)}  && \\
e^{ \frac{y}{c}}  &= \frac{x}{d}e^{ \frac{a}{c}(x+b)}  && \\
de^{ \frac{y}{c}}  &= xe^{ \frac{a}{c}x+\frac{ab}{c}}  && \\
de^{ \frac{y}{c}- \frac{ab}{c}}  &= xe^{ \frac{a}{c}x}  && \\
\frac{ad}{c}e^{ \frac{y}{c}- \frac{ab}{c}}  &= \frac{a}{c}xe^{ \frac{a}{c}x}  && \\
W( \frac{ad}{c}e^{ \frac{y}{c}- \frac{ab}{c}})  &= \frac{a}{c}x  && \\
\frac{W( \frac{ad}{c}e^{ \frac{y}{c}- \frac{ab}{c}})}{a}  &= x  && \\
\end{flalign*}
\begin{flalign*}
\intertext{Derivation \eqref{S4}, \eqref{C8}}
y  &= x^{a}+ \ln(x^{b}) && \\
e^{y}  &= x^{b}e^{x^{a}} && \\
e^{ \frac{y}{b}}  &= xe^{ \frac{x^{a}}{b}} && \\
e^{ \frac{ay}{b}}  &= x^{a}e^{ \frac{ax^{a}}{b}} && \\
e^{ \frac{ay}{b}}  &= \frac{a}{b}x^{a}e^{ \frac{ax^{a}}{b}} && \\
W(e^{ \frac{ay}{b}})  &= \frac{a}{b}x^{a} && \\
\frac{bW(e^{ \frac{ay}{b}})}{a}  &= x^{a} && \\
\left( \frac{bW(e^{ \frac{ay}{b}})}{a} \right)^{ \frac{1}{a}}  &= x && \\
\end{flalign*}
\end{minipage}
\newpage
\begin{minipage}[t]{0.5\textwidth}
\begin{flalign*}
\intertext{Derivation \eqref{S5}, \eqref{C9}}
y  &= \ln(x^{a})+W(x^{b}) && \\
e^{y}  &= e^{ \ln(x^{a})+W(x^{b})} && \\
e^{y}  &= x^{a}e^{W(x^{b})} && \\
e^{y}  &= x^{a}e^{W(x^{b})} && \\
u  &= x^{b} && \\
e^{y}  &= u^{ \frac{a}{b}}e^{W(u)} && \\
e^{y}  &= \left( W(u)e^{W(u)} \right)^{ a/b}e^{W(u)} && \\
e^{y}  &= W(u)^{ a/b}e^{( a/b+1)W(u)} && \\
e^{ \frac{by}{a}} &= W(u)e^{ \frac{a+b}{a}W(u)} && \\
\tfrac{a+b}{a}e^{ \frac{by}{a}} &= \tfrac{a+b}{a}W(u)e^{ \frac{a+b}{a}W(u)} && \\
W(\tfrac{a+b}{a}e^{ \frac{by}{a}}) &= \tfrac{a+b}{a}W(u) && \\
\frac{aW(\frac{a+b}{a}e^{ \frac{by}{a}})}{a+b} &= W(u) && \\
\frac{aW(\frac{a+b}{a}e^{ \frac{by}{a}})}{a+b}e^{\frac{aW(\frac{a+b}{a}e^{ \frac{by}{a}})}{a+b}} &= u && \\
\frac{aW(\frac{a+b}{a}e^{ \frac{by}{a}})}{a+b}e^{\frac{aW(\frac{a+b}{a}e^{ \frac{by}{a}})}{a+b}} &= x^{b} && \\
[\frac{aW(\frac{a+b}{a}e^{ \frac{by}{a}})}{a+b}]^{ \frac{1}{b}}e^{\frac{aW(\frac{a+b}{a}e^{ \frac{by}{a}})}{b(a+b)}} &= x && \\
\end{flalign*}
\end{minipage}
\begin{minipage}[t]{0.5\textwidth}
\begin{flalign*}
\intertext{Derivation \eqref{S6}}
y  &= W(ax)-W(bx) && \\
u &=ax, && \\
y &=W(u)-W( \tfrac{b}{a}u) && \\
e^{y} &= \tfrac{W( \tfrac{b}{a}u)}{bW(u)} && \\
e^{y}-\tfrac{a}{b} &= \frac{-a[W(u)-W( \tfrac{b}{a}u)]}{bW(u)}. && \\
e^{y}-\tfrac{a}{b} &= \frac{-a[W(u)-W( \tfrac{b}{a}u)]}{bW(u)} && \\
\frac{y}{e^{y}-\tfrac{a}{b}} &= \frac{W(u)-W( \tfrac{b}{a}u)}{ \frac{-a[W(u)-W( \tfrac{b}{a}u)]}{bW(u)}} && \\
&= bW(u). && \\
\frac{y}{be^{y}-a}e^{\frac{y}{be^{y}-a}} &= u = ax && \\
\frac{y}{abe^{y}-a^{2}} \exp \left(\frac{y}{be^{y}-a} \right) &= x
\end{flalign*}
\end{minipage}
\newpage
\begin{align*}
\intertext{Derivation \eqref{S7}}
y &= W(x)+W(1/x) && \\
y^{2}e^{y} &= [W(x)^{2}+2W(x)W(1/x)+W(1/x)^2]\frac{1}{W(x)W(1/x)} && \\
&= \frac{W(x)}{W(1/x)}+2+ \frac{W(1/x)}{W(1/x)} && \\
y^{2}e^{y}-2 &= \frac{W(x)}{W(1/x)}+ \frac{W(1/x)}{W(x)} && \\
\intertext{At this step, note that the right-hand side is in the form of $u+ \tfrac{1}{u}$ if $u= \frac{W(x)}{W(1/x)}$. The inverse operation of this can therefore be used to return $u$.} 
\tfrac{1}{2} \left[y^{2}e^{y}-2 \pm \sqrt{(y^{2}e^{y}-2)^{2}-4} \right] &= \left( \frac{W(x)}{W(1/x)} \right)^{ \pm 1}. && \\
\end{align*}
\begin{align*}
\intertext{There are now two possible inverse relations that will yield slightly different results. Take case 1 to be the positive branch of the square root, then,}
\tfrac{1}{2} \left[y^{2}e^{y}-2 + \sqrt{(y^{2}e^{y}-2)^{2}-4} \right] &= \frac{W(x)}{W(1/x)} && \\
e^{y} \tfrac{1}{2} \left[y^{2}e^{y}-2 + \sqrt{(y^{2}e^{y}-2)^{2}-4} \right] &= \frac{1}{W(x)W(1/x)} \frac{W(x)}{W(1/x)} && \\
&= \frac{1}{W(1/x)^{2}} && \\
2 e^{-y} \left[y^{2}e^{y}-2 + \sqrt{(y^{2}e^{y}-2)^{2}-4} \right]^{-1} &= W(1/x)^{2}  && \\
\sqrt{2}e^{-y/2} \left[y^{2}e^{y}-2 + \sqrt{(y^{2}e^{y}-2)^{2}-4} \right]^{-1/2} &= W(1/x)  && \\
\sqrt{2}e^{-y/2} \left[y^{2}e^{y}-2 + \sqrt{(y^{2}e^{y}-2)^{2}-4} \right]^{-1/2} e^{\sqrt{2}e^{-y/2} \left[y^{2}e^{y}-2 + \sqrt{(y^{2}e^{y}-2)^{2}-4} \right]^{-1/2}} &= 1/x  && \\
\frac{e^{y/2}}{ \sqrt{2}} \left[y^{2}e^{y}-2 + \sqrt{(y^{2}e^{y}-2)^{2}-4} \right]^{1/2} \exp \left(- \sqrt{2}e^{-y/2} \left[y^{2}e^{y}-2 + \sqrt{(y^{2}e^{y}-2)^{2}-4} \right]^{-1/2} \right) &= x  && \\
\intertext{Now to case 2 with the negative branch of the square root:}
\tfrac{1}{2} \left[y^{2}e^{y}-2 - \sqrt{(y^{2}e^{y}-2)^{2}-4} \right] &= \frac{W(1/x)}{W(x)} && \\
e^{y}\tfrac{1}{2} \left[y^{2}e^{y}-2 - \sqrt{(y^{2}e^{y}-2)^{2}-4} \right] &= \frac{1}{W(x)W(1/x)} \frac{W(1/x)}{W(x)} && \\
&= \frac{1}{W(x)^2} && \\
\frac{e^{-y}}{2} \left[y^{2}e^{y}-2 - \sqrt{(y^{2}e^{y}-2)^{2}-4} \right]^{-1} &= W(x)^{2} && \\
\frac{e^{-y/2}}{ \sqrt{2}} \left[y^{2}e^{y}-2 - \sqrt{(y^{2}e^{y}-2)^{2}-4} \right]^{-1/2} &= W(x) && \\
\frac{e^{-y/2}}{ \sqrt{2}} \left[y^{2}e^{y}-2 - \sqrt{(y^{2}e^{y}-2)^{2}-4} \right]^{-1/2} \exp \left( \frac{e^{-y/2}}{ \sqrt{2}} \left[y^{2}e^{y}-2 - \sqrt{(y^{2}e^{y}-2)^{2}-4} \right]^{-1/2} \right) &= x && \\
\end{align*}
\begin{align*}
\intertext{Derivation \eqref{S8}}
y &= 2W(ax)-W([bx]^{2}) && \\
u &= ax && \\
y &= W(u)-W([bu/a]^{2}) && \\
e^{y} &= \left[ \frac{u}{W(u)} \right]^{2} \frac{W([bu/a]^{2})}{   [ \frac{bu}{a}]^{2}} && \\
&= \frac{a^{2}W([bu/a]^2)}{b^{2}W(u)^{2}} && \\
\frac{b^{2}}{a^{2}}e^{y} &= \frac{W([bu/a]^2)}{W(u)^{2}} && \\
\frac{b^{2}}{a^{2}}ye^{y} &= [2W(u)-W([bu/a]^{2})]\frac{W([bu/a]^2)}{W(u)^{2}} && \\
\frac{b^{2}}{a^{2}}ye^{y}-1 &= \frac{2W(u)W([bu/a]^{2})-W([bu/a]^{2})^{2}}{W(u)^{2}}-\frac{W(u)^2}{W(u)^{2}} && \\
1-\frac{b^{2}}{a^{2}}ye^{y} &= \frac{[W(u)-W([bu/a]^{2}]^{2}}{W(u)^{2}} && \\
\pm \sqrt{1-\frac{b^{2}}{a^{2}}ye^{y}} &= \frac{W(u)-W([bu/a]]}{W(u)} && \\
\intertext{Case 1: }
\sqrt{1-\frac{b^{2}}{a^{2}}ye^{y}} &= \frac{W(u)-W([bu/a]^2)}{W(u)} && \\
\sqrt{1-\frac{b^{2}}{a^{2}}ye^{y}}-1 &= \frac{-W([bu/a]^2)}{W(u)} && \\
-\frac{a^{2}}{b^{2}}e^{-y} \left(\sqrt{1-\frac{b^{2}}{a^{2}}ye^{y}}-1 \right) &=\frac{W(u)^{2}}{W([bu/a]^{2})} \frac{W([bu/a]^{2})}{W(u)} && \\ 
&= W(u) && \\
-\frac{a^{2}}{b^{2}}e^{-y} \left(\sqrt{1-\frac{b^{2}}{a^{2}}ye^{y}}-1 \right)e^{-\frac{a^{2}}{b^{2}}e^{-y} \left(\sqrt{1-\frac{b^{2}}{a^{2}}ye^{y}}-1 \right)} &= u && \\
-\frac{a}{b^{2}}e^{-y} \left(\sqrt{1-\frac{b^{2}}{a^{2}}ye^{y}}-1 \right)e^{-\frac{a^{2}}{b^{2}}e^{-y} \left(\sqrt{1-\frac{b^{2}}{a^{2}}ye^{y}}-1 \right)} &= x && \\
\intertext{Case 2:}
\sqrt{1-\frac{b^{2}}{a^{2}}ye^{y}} &= \frac{W([bu/a]^2)-W(u)}{W(u)} && \\
\sqrt{1-\frac{b^{2}}{a^{2}}ye^{y}}+1 &= \frac{W([bu/a]^2)}{W(u)} && \\
\frac{a^{2}}{b^{2}}e^{-y} \left( \sqrt{1-\frac{b^{2}}{a^{2}}ye^{y}}+1 \right)e^{ \frac{a^{2}}{b^{2}}e^{-y} \left( \sqrt{1-\frac{b^{2}}{a^{2}}ye^{y}}+1 \right)} &= u && \\
\frac{a}{b^{2}}e^{-y} \left( \sqrt{1-\frac{b^{2}}{a^{2}}ye^{y}}+1 \right)e^{ \frac{a^{2}}{b^{2}}e^{-y} \left( \sqrt{1-\frac{b^{2}}{a^{2}}ye^{y}}+1 \right)} &= x && \\
\end{align*}
\begin{align*}
\intertext{Derivation \eqref{S9}}
y &= \sum_{k=0}^{n}a_{k}x \ln(b_{k}x) = a_{0}x \ln(b_{0}x)+a_{1}x \ln(b_{1}x)+...+a_{n}x \ln(b_{n}x) 
&& \\
&= x \sum_{k=0}^{n} \ln(b_{k}^{a_{k}}x^{a_{k}}) = x[a_{0} \ln(b_{0}x)+a_{1} \ln(b_{1}x)+...+a_{n} \ln(b_{n}x)] 
&&\\ 
e^{y} &= \exp \left(x \sum_{k=0}^{n} \ln(b_{k}^{a_{k}}x^{a_{k}}) \right) = \exp \left( x[ \ln(b_{0}^{a_{0}}x^{a_{0}})+ \ln(b_{1}^{a_{1}}x^{a_{1}})+...+ \ln(b_{n}^{a_{n}}x^{a_{n}})] \right) 
&&\\ 
&= \prod_{k=0}^{n}[b_{k}x]^{a^{k}x} = [b_{0}^{a_{0}}x^{a_{0}} \cdot b_{1}^{a_{1}}x^{a_{1}} \cdot... \cdot b_{n}^{a_{n}}x^{a_{n}}]^{x}
&& \\
&= \left[ \left( \prod_{k=0}^{n}b_{k}^{a_{k}} \right) x^{ \sum_{k=0}^{n}a_{k}} \right]^{x} = [b_{0}^{a_{0}}x^{a_{0}} \cdot b_{1}^{a_{1}}x^{a_{1}} \cdot... \cdot b_{n}^{a_{n}}x^{a_{n}}]^{x}
&& \\
u &= \sum_{k=0}^{n}a_{k}, && \\
v &= \prod_{k=0}^{n}b_{k}^{a_{k}}, 
&& \\
e^{y} &= [vx^{u}]^{x}
&& \\
e^{y/u} &= [v^{1/u}x]^{x}
&& \\
\exp \left( \frac{v^{1/u}y}{u} \right) &= [v^{1/u}x]^{v^{1/u}x}
&& \\
\exp \left( W \left( \ln \left( \exp \left( \frac{v^{1/u}y}{u} \right) \right) \right) \right) &= v^{1/u}x
&& \\
\frac{ \exp \left( W \left( \frac{v^{1/u}y}{u} \right) \right)}{v^{1/u}} &= x
&& \\
\frac{ \exp \left( W \left( \frac{ \left[ \prod_{k=0}^{n}b_{k}^{a_{k}} \right]^{1/ \sum_{k=0}^{n}a_{k}}y}{ \sum_{k=0}^{n}a_{k}} \right) \right)}{\left[ \prod_{k=0}^{n}b_{k}^{a_{k}} \right]^{1/ \sum_{k=1}^{n}a_{k}}} &= x
&& \\
\end{align*}
\newpage
\begin{minipage}[t]{0.5\textwidth}
\begin{align*}
\intertext{Derivation \eqref{C3}, \eqref{C7}}
y  &= (ax^{b})^{cx^d} && \\
y^{ \frac{1}{b}}  &= (a^{ \frac{1}{b}}x)^{cx^d} && \\
y^{ \frac{1}{bc}}  &= (a^{ \frac{1}{b}}x)^{x^d} && \\
y^{ \frac{d}{bc}}  &= (a^{ \frac{d}{b}}x^{d})^{x^d} && \\
y^{ \frac{a^{ \frac{d}{b}}d}{bc}}  &= (a^{ \frac{d}{b}}x^{d})^{a^{ \frac{d}{b}}x^d} && \\
sr_2 \left(y^{ \frac{a^{ \frac{d}{b}}d}{bc}} \right)  &= a^{ \frac{d}{b}}x^{d} && \\
e^{W( \frac{a^{ \frac{d}{b}}d}{bc} \ln(y))}  &= (a^{ \frac{1}{b}}x)^{d} && \\
e^{W( \frac{a^{ \frac{d}{b}}d}{bc} \ln(y))/d}  &= a^{ \frac{1}{b}}x && \\
\frac{e^{W( \frac{a^{ \frac{d}{b}}d}{bc} \ln(y))/d}}{a^{ \frac{1}{b}}}  &= x && \\
\end{align*}
\end{minipage}
\begin{minipage}[t]{0.5\textwidth}
\begin{align*}
\intertext{Derivation \eqref{C10}}
y &= a \tsr(bx)^{cx} && \\
\frac{y}{a} &= \tsr(bx)^{cx} && \\
\left[ \frac{y}{a} \right]^{b/c} &= \tsr(bx)^{bx} && \\
&= \tsr(bx)^{ \tsr(bx)^{ \tsr(bx)}} && \\
\tcr \left( \left[ \frac{y}{a} \right]^{b/c} \right) &= \tsr(bx) && \\
^{2} \left( \tcr \left( \left[ \frac{y}{a} \right]^{b/c} \right) \right) &= bx && \\
\frac{^{2} \left( \tcr \left( \left[ \frac{y}{a} \right]^{b/c} \right) \right)}{b} &= x && \\
\end{align*}
\end{minipage}
\begin{align*}
\intertext{Derivation \eqref{C4}}
y  &= a(bx^{ \ln(c)})^{ \ln(x)} && \\
y  &= e^{ \ln(a)}[e^{ \ln(b)}e^{ \ln(c) \ln(x)}]^{ \ln(x)} && \\
y  &= e^{ \ln(a)}[e^{ \ln(b) \ln(x)}e^{ \ln(c) \ln^{2}(x)}] && \\
y  &= e^{ \ln(c) \ln^{2}(x)+ \ln(b) \ln(x)+ \ln(a)} && \\
\ln(y)  &= \ln(c) \ln^{2}(x)+ \ln(b) \ln(x)+ \ln(a) && \\
0  &= \ln(c) \ln^{2}(x)+ \ln(b) \ln(x)+ \ln(a)- \ln(y) && \\
\frac{- \ln(b) \pm \sqrt{ \ln^{2}(b)-4 \ln(c)[ \ln(y)- \ln(a)]}}{2 \ln(c)}  &= \ln(x) && \\
\exp \left( \frac{- \ln(b) \pm \sqrt{ \ln^{2}(b)-4 \ln(c)[ \ln(y)- \ln(a)]}}{2 \ln(c)} \right)  &= x && \\
\end{align*}
\newpage
\begin{align*}
\intertext{Derivation \eqref{C6}}
y  &= ax^{b \tsr(x^{c})} && \\
\tfrac{y}{a}  &= x^{b \tsr(x^{c})} && \\
\left (\tfrac{y}{a} \right)^{ c/b}  &= x^{c \tsr(x^{c})} && \\
\left( \tfrac{y}{a} \right)^{ c/b}  &= e^{c \ln(x)\tsr(x^{c})} && \\
\left(\tfrac{y}{a} \right)^{ c/b}  &= e^{c \ln(x)e^{W(c \ln(x))}} && \\
\ln \left( \left( \tfrac{y}{a} \right)^{ c/b} \right)  &= c \ln(x)e^{W(c \ln(x))} && \\
\ln(( \tfrac{y}{a})^{ c/b})  &= W(c \ln(x))e^{W(c \ln(x))}e^{W(c \ln(x))} && \\
2 \ln(( \tfrac{y}{a})^{ c/b})  &= 2W(c \ln(x))e^{2W(c \ln(x))} && \\
W(2 \ln(( \tfrac{y}{a})^{ c/b}))  &= 2W(c \ln(x)) && \\
\frac{1}{2}W(2 \ln(( \tfrac{y}{a})^{ c/b}))  &= W(c \ln(x)) && \\
\frac{1}{2}W(2 \ln(( \tfrac{y}{a})^{ c/b}))e^{\frac{1}{2}W(2 \ln(( \tfrac{y}{a})^{ c/b}))}  &= c \ln(x) && \\
\frac{1}{2c}W(2 \ln(( \tfrac{y}{a})^{ c/b}))^{ \frac{1}{2}}(2 \ln(( \tfrac{y}{a})^{ \frac{c}{a}}))^{ \frac{1}{2}} &= \ln(x) && \\
\exp \left( \frac{1}{ \sqrt{2}c}(W(2 \ln(( \frac{y}{a})^{ c/b})) \ln(( \tfrac{y}{a})^{ \frac{c}{a}})^{1/2} \right) &= x && \\
\end{align*}
\newpage
\begin{align*}
\intertext{Derivation \eqref{Exa1}}
x+f+a &= ae^{b[y-xc]}+\frac{y-xc}{d} && \\
db[x+f+a] &= abde^{b[y-xc]}+b[y-xc] && \\
e^{db[x+f+a]} &= e^{b[y-xc]}e^{abde^{b[y-xc]}} && \\
abde^{db[x+f+a]} &= abde^{b[y-xc]}e^{abde^{b[y-xc]}} && \\
W \left(abde^{db[x+f+a]} \right) &= abde^{b[y-cx]} && \\
W \left(abde^{db[x+f+a]} \right) &= abde^{b[y-cx]} && \\
\ln(W \left(abde^{db[x+f+a]} \right)) &= \ln(abd)+b[y-cx] && \\
\ln \left(abde^{db[x+f+a]} \right)-W \left(abde^{db[x+f+a]} \right) &= && \\
\ln(abd)+db[x+f+a]-W \left(abde^{db[x+f+a]} \right) &= && \\
db[x+f+a]-W \left(abde^{db[x+f+a]} \right) &= b[y-cx] && \\
d[x+f+a]- \tfrac{1}{b}W \left(abde^{db[x+f+a]} \right) &= y-cx && \\
x[d+c]+d[f+a]- \tfrac{1}{b}W \left(abde^{db[x+f+a]} \right) &= y && \\
I[R_{sh}+R_{s}]+R_{sh}[I_{ph}+I_{0}]-\tfrac{n K_{b}T}{q}W \left(\tfrac{qI_{0}R_{sh}}{n K_{b}T} \exp \left( \frac{ qR_{sh}[I+I_{ph}+I_{0}]}{(n K_{b}T)} \right) \right) &= V && \\
\end{align*}
\begin{align*}
\intertext{ Derivation \eqref{Exa2}} D(t) &= \frac{b}{a} \left[ 1+W(-e^{-1- a^{2}t/b^{2}}) \right] && \\
\frac{aD}{b}-1 &= W(-e^{-1-a^{2}t/b^{2}}) && \\
[ \tfrac{a}{b}D-1]e^{[ \tfrac{a}{b}D-1]} & =-e^{-1-a^{2}t/b^{2}} && \\
\ln \left(-[ \tfrac{a}{b}D-1]e^{[ \tfrac{a}{b}D-1]} \right) &= -1-\frac{a^{2}t}{b^{2}} && \\
-[ \tfrac{a}{b}D-1]+\ln( \tfrac{a}{b}D-1)+1 &= -\frac{a^{2}t}{b^{2}} && \\
\tfrac{b^{2}}{a^2}[ \tfrac{a}{b}D-1)]+ \tfrac{b^2}{a^2} \ln( \tfrac{a}{b}D-1)+ \tfrac{b^{2}}{a^{2}} &= t
\end{align*}
\newpage
\begin{align*}
\intertext{Derivation \eqref{CTP1}}
y &= Ae^{Bx} && \\
\rightarrow r \sin( \theta+ \phi) &= Ae^{Br \cos( \theta+\phi)} && \\
r \sin( \theta) \cos(\phi)+r \cos( \theta) \sin(\phi) &= Ae^{B[r \cos( \theta) \cos(\phi)-r \sin( \theta) \sin(\phi)]} && \\
u &= r \cos(\theta), &&\\
a &= \cos(\phi), &&\\
v &= r \sin(\theta), &&\\
b &= \sin(\phi), &&\\
av+bu &= Ae^{Bau-Bbv} && \\
bu &= Ae^{Bau-Bbv}-av && \\
\tfrac{b}{a}u &= \tfrac{A}{a}e^{Bau-Bbv}-v && \\
\tfrac{Bb^2}{a}u &= \tfrac{ABb}{a}e^{Bau-Bbv}-Bbv && \\
\tfrac{Bb^2}{a}u &= \tfrac{ABb}{a}e^{Bau-Bbv}-Bbv && \\
\tfrac{Bb^2}{a}u+Bau &= \tfrac{ABb}{a}e^{Bau-Bbv}+Bau-Bbv && \\
e^{ B\frac{a^2+b^2}{a}u} &= e^{ \frac{ABb}{a}e^{Bau-Bbv}+Bau-Bbv} && \\
\tfrac{ABb}{a}e^{ B \frac{a^2+b^2}{a}u} &= \tfrac{ABb}{a}e^{Bau-Bbv}e^{ \frac{ABb}{a}e^{Bau-Bbv}} && \\
W( \tfrac{ABb}{a}e^{ B \frac{a^2b^2}{a}u}) &= \tfrac{ABb}{a}e^{Bau-Bbv} && \\
\tfrac{a}{ABb}W( \tfrac{ABb}{a}e^{ B \frac{a^2+b^2}{a}u}) &= e^{Bau-Bbv} && \\
\ln( \tfrac{a}{ABb}W( \tfrac{ABb}{a}e^{ B\frac{a^2+b^2}{a}u})) &= Bau-Bbv && \\
\ln( \tfrac{a}{ABb}W( \tfrac{ABb}{a}e^{ B\frac{a^2+b^2}{a}u}))-Bau &= Bbv && \\
\tfrac{1}{Bb} \ln( \tfrac{a}{ABb}W( \tfrac{ABb}{a}e^{ B\frac{a^2+b^2}{a}u}))- \frac{a}{b}u &= v && \\
\frac{1}{B \sin(\phi)} \ln \left( \frac{ \cos(\phi)}{AB \sin(\phi)}W( \frac{AB \sin(\phi)}{ \cos(\phi)}e^{ B\frac{ \cos(\phi)^2+ \sin(\phi)^2}{ \cos(\phi)}r \cos( \theta)}) \right)- \frac{ \cos(\phi)}{ \sin(\phi)}r \cos( \theta) &= r \sin( \theta) && \\
\frac{ \csc(\phi)}{B} \ln \left( \frac{ \cot(\phi)}{AB}W(AB \tan(\phi)e^{B \sec(\phi)x}) \right)- \cot(\phi)x &= y && \\
\end{align*}
\newpage
\label{End Proofs} \begin{align*}
\intertext{Derivation \eqref{CTP2}}
y &= \frac{1}{B} \ln( \frac{x}{A}) && \\
r \sin( \theta) &= \frac{1}{B} \ln( \frac{r \cos( \theta)}{A}) && \\
\rightarrow r \sin( \theta+\phi) &= \frac{1}{B} \ln( \frac{r \cos( \theta+\phi)}{A}) && \\
r \sin( \theta) \cos(\phi)+ r \cos( \theta) \sin(\phi) &= \frac{1}{B} \ln( \frac{r \cos( \theta) \cos(\phi)- r \sin( \theta) \sin(\phi)}{A}) && \\
u = r \cos(\theta), a = \cos(\phi), v = r \sin(\theta), b = \sin(\phi) &&\\
v a+ u b &= \frac{1}{B} \ln( \frac{ua-vb}{A}) && \\
e^{av+ bu} &= ( \frac{au-bv}{A})^{ \frac{1}{B}} && \\
e^{Bav+ Bbu} &= ( \frac{au-bv}{A}) && \\
Ae^{Bav+ Bbu} &= au-bv && \\
\frac{A}{b}e^{Bav+ Bbu} &= \frac{a}{b}u-v && \\
\frac{A}{b}e^{Bav+ Bbu}+v &= \frac{a}{b}u && \\
\frac{ABa}{b}e^{Bav+ Bbu}+Bav &= \frac{Ba^{2}}{b}u && \\
\frac{ABa}{b}e^{Bav+ Bbu}+Bav &= \frac{Ba^{2}}{b}u+Bbu && \\
\frac{ABa}{b}e^{Bav+ Bbu}+Bav+Bbu &= \frac{Ba^{2}+Bb^{2}}{b}u && \\
e^{\frac{ABa}{b}e^{Bav+ Bbu}+Bav+Bbu} &= e^{\frac{Ba^{2}+Bb^{2}}{b}u} && \\
\frac{ABa}{b}e^{Bav+Bbu}e^{\frac{ABa}{b}e^{Bav+ Bbu}} &= \frac{ABa}{b}e^{\frac{Ba^{2}+Bb^{2}}{b}u} && \\
\frac{ABa}{b}e^{Bav+Bbu} &= W( \frac{ABa}{b}e^{\frac{Ba^{2}+Bb^{2}}{b}u}) && \\
e^{Bav+Bbu} &= \frac{b}{ABa} W( \frac{ABa}{b}e^{\frac{Ba^{2}+Bb^{2}}{b}u}) && \\
Bav+Bbu &= \ln( \frac{b}{ABa} W( \frac{ABa}{b}e^{\frac{Ba^{2}+Bb^{2}}{b}u})) && \\
Bav &= \ln( \frac{b}{ABa} W( \frac{ABa}{b}e^{\frac{Ba^{2}+Bb^{2}}{b}u}))-Bbu && \\
v &= \frac{1}{Ba} \ln( \frac{b}{ABa} W( \frac{ABa}{b}e^{\frac{Ba^{2}+Bb^{2}}{b}u}))- \frac{b}{a}u && \\
r \sin( \theta) &= \frac{ \sec(\phi)}{B} \ln( \frac{ \tan(\phi)}{AB}W(AB \cot(\phi) e^{B \csc(\phi)r \cos( \theta)}))- \tan(\phi)r \cos( \theta) && \\
y &= \frac{ \sec(\phi)}{B} \ln \left( \frac{ \tan(\phi)}{AB}W(AB \cot(\phi) e^{B \csc(\phi)x}) \right)- \tan(\phi)x && \\
\end{align*}
\newpage
\begingroup 
\makeatletter
\renewcommand{\@seccntformat}[1]{}
\makeatother 
\noindent \textbf{\Large Acknowledgments} \hfill \\
The author is grateful to the people who supported the author's passion for this paper. The author would like to publicly thank: David Dollhopf, engineer, for his unparalleled enthusiasm and personal support; David Jeffrey, professor of applied mathematics for his exceptional patience and personal correspondence; Suzanne Woll, chemical engineer, for her time in providing feedback and guidance; Alex Hanhart, professor of mathematics for his continuous guidance in the author's pursuit of mathematics. 
\endgroup
\invisiblesection{References} 

\end{document}